\documentclass[preprint,12pt,3p]{elsarticle}

\usepackage{amssymb}
\usepackage{amsmath,amsthm}
\usepackage{mathrsfs}
\usepackage{hyperref}
\usepackage{cleveref}
\usepackage[all,cmtip]{xy}
\usepackage{epsf, graphicx}
\usepackage{latexsym,amsfonts,amsbsy,amssymb}
\usepackage{geometry}
\usepackage{titletoc}
\usepackage{rotating}
\usepackage{algorithm}
\usepackage{algorithmic}
\usepackage{amscd}

\makeatletter

\@addtoreset{equation}{section} \makeatother
\newtheorem{theorem}{Theorem}[section]
\newtheorem{lemma}[theorem]{Lemma}

\newtheorem{corollary}[theorem]{Corollary}
\newtheorem{remark}[theorem]{Remark}

\newtheorem{proposition}[theorem]{Proposition}

\def\Ind{{\rm Ind}}

\def\Br{{\rm Br}}

\def\Aut{{\rm Aut}}

\def\Hom{{\rm Hom}}
\def\End{{\rm End}}

\def\G{\mathcal{G}}

\def\F{\mathcal{F}}
\def\G{\mathcal{G}}
\def\FF{\mathbb{F}}

\def\P{\mathcal{P}}

\def\Proj{{\rm Proj}}
\def\FFTD{\mathbb{F}T^\Delta}

\makeatletter
\def\ps@pprintTitle{%
\let\@oddhead\@empty
\let\@evenhead\@empty
\def\@oddfoot{\reset@font\hfil\thepage\hfil}
\let\@evenfoot\@oddfoot
}
\makeatother

\begin{document}

\begin{frontmatter}

\title{Stable equivalences with endopermutation source, slash functors, and functorial equivalences}

\author[label1]{Xin Huang}
\fntext[label1]{The author acknowledges support by China Postdoctoral Science Foundation (2023T160290).}

\address{SICM, Southern University of Science and Technology, Shenzhen 518055, China}
\ead{xinhuang@mails.ccnu.edu.cn}

\begin{abstract}
We show that a bimodule of two block algebras of finite groups which has an endopermutation module as a source and which induces a stable equivalence of Morita type gives rise, via slash functors, to a family of bimodules of local block algebras with endopermutation source which induce Morita equivalences. As an application, we show that Morita (resp. stable) equivalences with endopermutation source imply functorial (resp. stable functorial) equivalences defined by Bouc and Y{\i}lmaz.
\end{abstract}


\end{frontmatter}


\section{Introduction}\label{s1}

Throughout this paper, $p$ is a prime and $k$ is an algebraically closed  field of characteristic $p$.
Following Brou\'{e} (\cite[\S5.A]{Broue94}), for two symmetric $k$-algebras $A$, $B$ and an $(A,B)$-bimodule $M$, we say that $M$ induces a {\it stable equivalence of Morita type} between $A$ and $B$ if $M$ is finitely generated projective as a left $A$-module and as a right $B$-module, and if $M\otimes_BM^*\cong A\oplus W$ for some projective $A\otimes_k A^{\rm op}$-module $W$ and $M^*\otimes_A M\cong B\oplus W'$ for some projective $B\otimes_k B^{\rm op}$-module $W'$, where $M^*:={\rm Hom}_k(M,k)$.

For $P$, $Q$ finite $p$-groups, $\mathcal{F}$ a fusion system on $P$ and $\varphi:P\to Q$ a group isomorphism, denote by ${}^\varphi\F$ the fusion system on $Q$ induced by $\F$ via the isomorphism $\varphi$.
We set $\Delta \varphi:=\{(u,\varphi(u))|u\in P\}$ and whenever useful, we regard a $k\Delta \varphi$-module $V$ as a $kP$-module and vice versa via the isomorphism $P\cong \Delta \varphi$ sending $u\in P$ to $(u,\varphi(u))$.

\begin{theorem}[cf. {{\cite[7.6]{Puig1999}}, see \cite[Theorem 9.11.2]{Linckelmann}}]\label{Puig1}
Let $G$ and $H$ be finite groups, $b$ a block  of $kG$ and $c$ a block  of $kH$. Let $P$ be a defect group of $b$ and $R$ a defect group of $c$.
Let $i\in$ $(k Gb)^P$ and
$j\in$ $(k Hc)^R$ be source idempotents. Denote by $\F$ the fusion
system on $P$ of $b$ determined by $i$, and denote by $\mathcal{G}$ the fusion
system on $R$ determined by $j$. Let $M$ be an indecomposable
$(k Gb, k Hc)$-bimodule inducing a stable equivalence of Morita type between $k Gb$ and $k Hc$.
Suppose that as a $k(G\times H)$-module, $M$ has an endopermutation module as a source for some vertex.
Then there is an isomorphism $\varphi : P\to R$ and an indecomposable
endopermutation $\Delta\varphi$-module $V$ such that $M$ is isomorphic
to a direct summand of
$$k G i\otimes_{k P} \Ind_{\Delta\varphi}^{P\times R}(V) \otimes_{k R} jk H$$
as a $(k Gb,k Hc)$-bimodule. For any such $\varphi$ and
$V$, the following hold.

\begin{enumerate} [\rm (a)]

\item
The group $\Delta\varphi$ is a vertex of $M$, the group $\Delta \varphi^{-1}$ is a vertex of its dual $M^*$, the module $V$ is a $k\Delta\varphi$-source of $M$ and its dual $V^*$ is a $k\Delta\varphi^{-1}$-source of $M^*$.

\item
We have ${^\varphi{\F}}=$ $\mathcal{G}$, and when regarded as a $k P$-module, the endopermutation module $V$ is $\F$-stable.
\end{enumerate}
\end{theorem}

Keep the notation and assumptions of Theorem \ref{Puig1}. We may identify $P$ and $R$ via $\varphi$, so now $P=R$, $\F=\mathcal{G}$, and $\Delta \varphi=\Delta P:=\{(u,u)|u\in P\}$. (This simplifies notation, but one
could as well go on without making those identifications.) For any non-trivial subgroup $Q$ of $P$, denote by $e_Q$ and $f_Q$ the unique blocks of $kC_G(Q)$ and $kC_H(Q)$, respectively, satisfying ${\rm br}_Q^{kG}(i)e_Q\neq 0$ and ${\rm br}_Q^{kH}(j)f_Q\neq 0$. Let $N_G(Q,e_Q):=\{g\in N_G(Q)~|~ge_Qg^{-1}=e_Q\}$.  Since $\F=\mathcal{G}$, we have $N_G(Q,e_Q)/C_G(Q)\cong N_H(Q,f_Q)/C_H(Q)$. Let $G_Q$ be an arbitrary subgroup of $N_G(Q,e_Q)$ containing $C_G(Q)$, then there is a subgroup $H_Q$ of $N_H(Q,f_Q)$ such that $H_Q/C_H(Q)$ is the image of $G_Q/C_G(Q)$ under the above isomorphism. By \cite[Theorem 6.2.6 (iii)]{Linckelmann}, $e_Q$ (resp. $f_Q$) remains a block of $kG_Q$ (resp. $kH_Q$).

In \cite{Biland}, Biland introduced the notion of Brauer-friendly modules, a generalisation of endo-$p$-permutation modules. Biland showed that Dade's slash construction for endopermutation modules (see \cite[Section 4 and 5]{Dade}) can be extended for applying to Brauer-friendly modules and can also be turned into a functor on a category of compatible Brauer-friendly modules, sharing most of the very useful properties that are satisfied by the Brauer functor over the category of $p$-permutation modules. Biland expected in \cite{Biland} that the slash functors may have important consequences in local representation theory. For instance, it would make it easier to deal with the local properties of derived equivalences between blocks algebras, even though this equivalence is not splendid.

The purpose of this paper is to use the slashed modules constructed by Biland to give quick and explicit bimodules which induce the Morita equivalences between local blocks. This confirms Biland's expectation. Consider the block $b\otimes c^\circ$ of $k(G\times H)$, where $c^\circ$ is the image of $c$ in $kH$ under the anti-automorphism of $kH$ sending $h\in H$ to $h^{-1}$. By Proposition \ref{brauerfriendly} below, $M$ is a Brauer-friendly $k(G\times H)(b\otimes c^{\circ})$-module. Denote by $\Delta Q$ the subgroup $\{(u,u)|u\in Q\}$ of $\Delta P$. Then the pair $(\Delta Q, e_Q\otimes f_Q^\circ)$ is a $(b\otimes c^\circ)$-subpair. Let $M_Q$ be a $(\Delta Q, e_Q\otimes f_Q^\circ)$-slashed module attached to $M$ over the group $N_{G\times H}(\Delta Q,e_Q\otimes f_Q^\circ)$ (see Lemma \ref{slash} below for the definition).
So $M_Q$ is a Brauer-friendly $kN_{G\times H}(\Delta Q, e_Q\otimes f_Q^\circ)$-module.
We show that the module $M_Q$ induces Morita equivalences between local blocks.

\begin{theorem}\label{main}
Keep the notation and assumptions above. The following hold.
\begin{enumerate} [\rm (a)]

\item The $(kC_G(Q)e_Q, kC_H(Q)f_Q)$-bimodule ${\rm Res}_{C_G(Q)\times C_H(Q)}^{N_{G\times H}(\Delta Q, e_Q\otimes f_Q^\circ)}(M_Q)$ induces a Morita equivalence between $kC_G(Q)e_Q$ and $kC_H(Q)f_Q$ and has an endopermutation module as a source.

\item Let $T_Q:=N_{G\times H}(\Delta Q, e_Q\otimes f_Q^\circ)\cap(G_Q\times H_Q)$. The $(kG_Qe_Q,kH_Qf_Q)$-bimodule ${\rm Ind}_{T_Q}^{G_Q\times H_Q}(M_Q)$ induces a Morita equivalence between $kG_Qe_Q$ and $kH_Qf_Q$ and has an endopermutation module as a source.

\end{enumerate}
\end{theorem}

Clearly the statement (b) contains (a). Actually, (b) can be easily deduced from (a) as well. The existence of a Morita equivalence between $kC_G(Q)e_Q$ and $kC_H(Q)f_Q$ with endopermutation source is a theorem of Puig (see \cite[7.7.4]{Puig1999}). Theorem \ref{main} (a) gives an explicit bimodule inducing such a Morita equivalence. Theorem \ref{main} (b) generalises results of Puig-Zhou (\cite[Theorem 1.4]{PZ}) and Hu (\cite[Theorem 1.1]{Hu}). By \cite[Theorem 40.4 (a)]{Thevenaz}, the group $N_G(Q_\delta)$ in \cite[Theorem 1.4]{PZ} is a subgroup of $N_G(Q,e_Q)$ in Theorem \ref{main}. The results of \cite{PZ} and \cite{Hu} start from a Morita equivalence while Theorem \ref{main} starts from a stable equivalence. Our proof of Theorem \ref{main} does not depend on either \cite[7.7.4]{Puig1999} or \cite[Theorem 1.4]{PZ}. Besides, Theorem \ref{main} (a) and \cite[Theorem 1.2]{Lin15} are the converses of each other. 

Since Theorem \ref{main} gives explicit local Morita equivalences in terms of slashed modules and since slash functors share most of the very useful properties of Brauer functors, we can obtain more information about this family of local Morita equivalences, and hence Theorem \ref{main} might have more applications. For example, since Theorem \ref{main} is an analogue of \cite[Theorem 1.3]{BP} for stable equivalences with endopermutation source instead of $p$-permutation equivalences, and since \cite[Theorem 1.3]{BP} played an important role in the proof of the main result of \cite{Boltje}, Theorem \ref{main} might be useful in determining the Brou\'{e} invariant of a Morita equivalence with endopermutation source. For another example, we can use Theorem \ref{main} (a) to prove the following result.

\begin{corollary}\label{functorialequ}
Let $\mathbb{F}$ be an algebraically closed filed of characteristic 0. Let $G$ and $H$ be finite groups and let $b$ and $c$ be blocks of $kG$ and $kH$ respectively. If there is a Morita (resp. stable) equivalence between $kGb$ and $kHc$ induced by a $(kGb,kHc)$-bimodule with endopermutation source, then $(G,b)$ and $(H,c)$ are functorially (resp. stably functorially) equivalent over $\mathbb{F}$ in the sense of Bouc and Y{\i}lmaz \cite{BY} (resp. \cite{BY2}).
\end{corollary}

Combined with \cite[Theorem 8.11.5 (i)]{Linckelmann} and \cite[Theorem 10.5.1]{Linckelmann}, Corollary \ref{functorialequ} recovers the following result of Bouc and Y{\i}lmaz.

\begin{corollary}[{\cite[Theorem 10.5 (v)]{BY} and \cite[Theorem 1.3]{BY2}}] Let $\mathbb{F}$ be an algebraically closed filed of characteristic 0.
\begin{enumerate} [\rm (a)]
	\item Let $G$ be a finite group and $b$ a nilpotent block of $kG$ with a defect group $P$. Then $(G,b)$ and $(P,1)$ are functorially equivalent over $\FF$.

\item Let $G$ be a finite group, $b$ a block of $kG$, $P$ a defect group of $b$, $i\in (kGb)^P$ a source idempotent of $b$ and $e_P$ the block of $kC_G(P)$ satisfying ${\rm br}_P^{kG}(i)e_P\neq 0$. Suppose that $P$ is non-trivial abelian and that $E:=N_G(P,e_P)/C_G(P)$ acts freely on $ P\setminus\{1\}$. Then $(G, b)$ and $(P\rtimes E,1)$ are stably functorially equivalent over $\FF$.
\end{enumerate}
\end{corollary}

Using a property of slash functors, we can obtain one more corollary of Theorem \ref{main}.

\begin{corollary}\label{coro1}
Keep the notation and assumptions of Theorem \ref{main}. The following hold.
\begin{enumerate} [\rm (a)]

\item ${\rm Res}_{C_G(P)\times C_H(P)}^{N_{G\times H}(\Delta P, e_P\otimes f_P^\circ)}(M_P)$ is a $p$-permutation module and induces a Morita equivalence between $kC_G(P)e_P$ and $kC_H(P)f_P$.

\item ${\rm Ind}_{T_P}^{G_P\times H_P}(M_P)$ is a $p$-permutation module and induces a Morita equivalence between $kG_Pe_P$ and $kH_Pf_P$.

\item ${\rm Ind}_{N_{G\times H}(\Delta P, e_P\otimes f_P^\circ)}^{N_G(P)\times N_H(P)}(M_P)$ is a $p$-permutation module and induces a Morita equivalence between the Brauer correspondent blocks of $b$ and $c$.
\end{enumerate}
\end{corollary}

As pointed out by Linckelmann (the paragraph below \cite[Theorem 1.2]{Lin15}), it is not clear whether $e_QMf_Q$ is an endopermutation $k\Delta Q$-module. To overcome this issue, we need to reformulate Theorem \ref{main} (a) at the level of almost source algebras.

\begin{theorem}\label{main2}
Let $A$ and $B$ be almost source algebras of blocks of finite group algebras having a common defect group $P$ and the same fusion system $\mathcal{F}$ on $P$. Let $V$ be an $\mathcal{F}$-stable indecomposable endopermutation $kP$-module with vertex $P$. Let $M$ be an indecomposable direct summand of the $(A,B)$-bimodule
$$A\otimes_{kP}{\rm Ind}_{\Delta P}^{P\times P}(V)\otimes_{kP}B.$$
Then for any non-trivial subgroup $Q$ of $P$, by Proposition \ref{lin154.1} below, there is a canonical $({\rm Br}_Q(A), \Br_Q(B))$-bimodule $M_Q$ satisfying ${\rm End}_k(M_Q)\cong {\rm Br}_{\Delta Q}({\rm End}_k(M))$.
Suppose that M and its dual $M^*$ induce a stable equivalence of Morita type between $A$ and $B$. Then for any non-trivial subgroup $Q$ of $P$, $M_Q$ induces a Morita equivalence between ${\rm Br}_Q(A)$ and $\Br_Q(B)$.
\end{theorem}

Theorem \ref{main2} and \cite[Theorem 1.1]{Lin15} are the converses of each other.
The terminology and required background information for the results above are collected in \S\ref{s2}.

\section{Preliminaries}\label{s2}

In this section we review some notation and background material. For an algebra $A$, we denote by $A^{\rm op}$ the opposite algebra of $A$. Unless specified otherwise, modules in the paper are left modules and are finite-dimensional over $k$. For a finite group $G$, we denote by $\Delta G$ the diagonal subgroup $\{(g,g)|g\in G\}$ of the direct product $G\times G$. For finite groups $G$ and $H$, a $(k G,kH)$-bimodule $M$ can be regarded as a $k(G\times H)$-module (and vice versa) via $(g,h)m=gmh^{-1}$, where $g\in G$, $h\in H$ and $m\in M$. If $M$ is indecomposable as a $(k G,k H)$-bimodule, then $M$ is indecomposable as a $k(G\times H)$-module, hence has a vertex (in $G\times H$) and a source.

Given two $k$-algebras $A$, $B$ and an $(A,B)$-bimodule $M$, ${\rm End}_k (M)$ is an $(A\otimes_k A^{\rm op},B\otimes_kB^{\rm op})$-bimodule: for any $a_1,a_2\in A$, $b_1,b_2\in B$, $\varphi\in {\rm End}_k(M)$ and $m\in M$,
$$((a_1\otimes a_2)\cdot\varphi\cdot(b_1\otimes b_2))(m)=a_1\varphi(a_2mb_2)b_1.$$

\subsection{The Brauer construction}

Let $G$ be a finite group. We can refer to \cite[\S10]{Thevenaz} or \cite[Definition 1.3.1]{Linckelmann} for the definition of a $G$-algebra. If $A$ is a $G$-algebra (resp. $kG$-module), we denote by $A^H$ the subalgebra (resp. submodule) of $H$-fixed points of $A$ for any subgroup $H$ of $G$.  For any two $p$-subgroups $Q\leq P$ of $G$, the {\it relative trace map} ${\rm Tr}_Q^P:A^Q\to A^P$, is defined by ${\rm Tr}_Q^P(a)=\sum_{x\in [P/Q]}{}^xa$, where $[P/Q]$ denotes a set of representatives of the left cosets of $Q$ in $P$. We denote by ${\rm Br}_P(A)$ the {\it Brauer quotient} of $A$, i.e., the $N_G(P)$-algebra (resp. $kN_G(P)$-module)
\[{\rm Br}_P(A)=A^P/\sum_{Q<P}{\rm Tr}_Q^P(A^Q).\]
We denote by ${\rm br}_P^A:A^P\to {\rm Br}_P(A)$ the canonical map, which is called the {\it Brauer homomorphism}.

If $A$ is a $G$-interior algebra (for instance, $A={\rm End}_k(M)$ for some $k G$-module $M$), then the Brauer quotient ${\rm Br}_P(A)$ has a natural structure of $C_G(P)$-interior algebra. For the group algebra $k G$ (considered as a $G$-interior algebra) and a $p$-subgroup $P$, the Brauer homomorphism can be identified with the $k$-algebra homomorphism
$(k G)^P\to kC_G(P)$, $\sum_{g\in G}\alpha_gg\mapsto \sum_{g\in C_G(P)}\alpha_gg$.
If $M$ is an $A$-module then $M$ can be viewed as a $kG$-module via the structure homomorphism $G\to A^\times$. The $A$-module structure on $M$ then induces a ${\rm Br}_P(A)$-module structure on the corresponding Brauer quotient ${\rm Br}_P(M)$. If $B$ is another $G$-interior algebra, and $M$ is an $(A,B)$-bimodule, then $M$ can be viewed as a $(kG,kG)$-bimodule, hence as a $k(G\times G)$-module. The $(A,B)$-bimodule structure on $M$ restricts to an $(A^P,B^P)$-bimodule structure on $M^{\Delta P}$ which in turn induces a $(\Br_P(A),\Br_P(B))$-bimodule structure on $\Br_{\Delta P}(M)$.

\begin{remark}\label{remark1}
{\rm Let $A$ be a $P$-interior algebra for some finite $p$-group $P$. On the one hand, taking the $P$-Brauer quotient of the $P$-interior algebra $A$, we obtain an algebra ${\rm Br}_P(A)$. We can regard ${\rm Br}_P(A)$ as a $(\Br_P(A),\Br_P(A))$-bimodule and denote it by $M_1$. On the other hand, considering the $(A,A)$-bimodule $A$, we can regard $A$ as a $k(P\times P)$-module. Taking the $\Delta P$-Brauer quotient of the $k(P\times P)$-bimodule $A$, we obtain a $(\Br_P(A),\Br_P(A))$-bimodule $M_2:={\rm Br}_{\Delta P}(A)$. The bimodules $M_1$ and $M_2$ are actually equal, although they are obtained by different ways. So we can use the notation ${\rm Br}_P(A)$ to stand for ${\rm Br}_{\Delta P}(A)$.}
\end{remark}

\subsection{Blocks and subpairs}\label{weights}

Let $G$ be a finite group. By a {\it block} of the group algebra $k G$, we mean a primitive idempotent $b$ of the center of $k G$, and $k Gb$ is called a {\it block algebra} of $k G$. A {\it defect group} of $b$ is a maximal $p$-subgroup $P$ of $G$ such that ${\rm br}_P^{k G}(b)\neq 0$.

A {\it subpair} of the group $G$ is a pair $(Q,e_Q)$, where $Q$ is a $p$-subgroup of $G$ and $e_Q$ is a block of $kC_G(Q)$. The idempotent $e_Q$ remains a block of $kH$ whenever $H$ is a subgroup of $G$ satisfying $C_G(Q)\leq H\leq N_G(Q,e_Q)$. The subpair $(Q,e_Q)$ is a {\it$b$-subpair} if $e_Q{\rm br}_Q^{kG}(b)\neq 0$.

Let $G$ and $H$ be finite groups. Denote by $-^\circ$ the $k$-algebra isomorphism $k H\cong(k H)^{\rm op}$ sending any $h\in H$ to $h^{-1}$. Let $b$ and $c$ be blocks of $k G$ and $k H$ respectively. Clearly $c^\circ$ is a block of $k H$. Then a $(k Gb,k Hc)$-bimodule $M$ can be regarded as a $k(G\times H)$-module belonging to the block $b\otimes c^\circ$ of $k(G\times H)$, and vice versa. Here, we identify $b\otimes c^\circ$ and its image under the $k$-algebra isomorphism $k G\otimes_k k H\cong k(G\times H)$ sending $g\otimes h$ to $(g,h)$ for any $g\in G$ and $h\in H$.

\subsection{Almost source algebras}

For the concept of fusion systems, we follow the conventions of \cite[\S8.1]{Linckelmann}. Almost source algebras of a block was introduced by Linckelmann in \cite[Definition 4.3]{Lin08}. Let $G$ be a finite group, $b$ a block of $kG$ and $P$ a defect group of $b$. An idempotent $i\in (kGb)^P$ is called an {\it almost source idempotent} if ${\rm br}_P^{kG}(i)\neq 0$ for every subgroup $Q$ of $P$, there is a unique block $e_Q$ of $kC_G(Q)$ such that ${\rm br}_Q^{kG}(i)\in kC_G(Q)e_Q$. The $P$-interior algebra $ikGi$ is then called an {\it almost source algebra} of the block $b$. By \cite[Proposition 4.1]{Lin08}, there is a canonical Morita equivalence between the block algebra $kGb$ and the almost source algebra $ikGi$ sending a $kGb$-module $M$ to the $ikGi$-module $iM$. The choice of an almost source idempotent $i\in (kGb)^P$ determines a fusion system $\F$ on $P$ such that for any subgroups $Q$ and $R$ of $P$, the set $\Hom_\F(Q,R)$ is the set of all group homomorphisms $\varphi:Q\to R$ for which there is an element $x\in G$ satisfying $\varphi(u)=xux^{-1}$ for all $u\in Q$ and satisfying $xe_Qx^{-1}=e_{xQx^{-1}}$ (see \cite[Remark 4.4]{Lin08} or \cite[\S 8.7]{Linckelmann}). Note that we use here the blanket assumption that $k$ is large enough.  Moreover, a subgroup $Q$ of $P$ is fully $\F$-centralised if and only if $C_P(Q)$ is a defect group of the block $e_Q$ of $kC_G(Q)$. Given a subgroup $Q$ of $P$, it is always possible to find a subgroup $R$ of $P$ such that $Q\cong R$ in $\F$ and such that $R$ is fully $\F$-centralised.
The following result of Linckelmann explains why we will need to work with fully centralised subgroups and almost source idempotent rather than source idempotents.

\begin{proposition}[{\cite[Proposition 4.5]{Lin08}}]\label{almost}
Let $G$ be a finite group, $b$ a block of $kG$, $P$ a defect group of $b$, and $i\in (kGb)^P$ an almost source idempotent of $b$ with the associated almost source algebra $A=ikGi$. If $Q$ is fully $\F$-centralised subgroup of $P$, then ${\rm br}_Q^{kG}(i)$ is an almost source idempotent of $kC_G(Q)e_Q$ with the associated almost source algebra $\Br_Q(A)$. In particular, $kC_G(Q)e_Q$ and $\Br_Q(A)$ are Morita equivalent.
\end{proposition}

\subsection{Bimodules with fusion-stable endopermutation source}

In this subsection we review some results on bimodules with fusion-stable endopermutation source proved by Linckelmann (\cite[\S4]{Lin15}). Let $P$ be a finite $p$-group, $V$ and $W$ be two $kP$-modules. Then there is a natural action of the group $P\times P$ on the $k$-module ${\rm Hom}_k(V,W)$ and in particular on the $k$-module ${\rm End}_k(V)$: for  $(g_1,g_2)\in P$, $f\in {\rm Hom}_k(V,W)$ and $v\in V$, $((g_1,g_2)\cdot f)(v)=g_2f(g_1^{-1}v)$. If ${\rm End}_k(V)$ admits a $\Delta P$-stable $k$-basis, then $V$ is called an {\it endopermutation} $kP$-module, as defined in \cite{Dade}.

Let $P$ be a finite $p$-group and $\F$ a fusion system on $P$. Let $Q$ be a subgroup of $P$ and $V$ an endopermutation $kQ$-module. Following Linckelmann (\cite[Definition 3.1]{Lin15}), we say that $V$ is {\it $\F$-stable} if for any subgroup $R$ of $Q$ and any morphism $\varphi:R\to Q$ in $\F$, the sets of isomorphism classes of indecomposable direct summands with vertex $R$ of the $kR$-modules ${\rm Res}_R^Q(V)$ and ${}_\varphi V$ are equal (including the possibility that both sets may be empty). This is equivalent to saying that ${\rm Res}_R^Q(V)\oplus {}_\varphi V$ is an endopermutation $kR$-module (see \cite[Corollary 6.12]{Dade}).

\begin{proposition}[{\cite[Proposition 3.2]{Lin15}}]\label{slash2}
Let $A$ be an almost source algebra of a block of a finite group $G$ over $k$ with a defect group $P$ and fusion system $\F$ on $P$. Let $Q$ be a subgroup of $P$ and let $V$ be an $\F$-stable endopermutation $kQ$-module having an indecomposable direct summand with vertex $Q$. Set $U:=A\otimes_{kQ}V$. The following hold.

\begin{enumerate} [\rm (i)]

\item As a $kQ$-module, $U$ is an endopermutation module, and $U$ has a direct summand isomorphic to $V$.

\item Let $R$ be a subgroup of $Q$ and $M$ a direct summand of $U$, then there exists a $\Br_R(A)$-module $M_Q$ such that we have an isomorphism $\Br_R(\End_k(M))\cong \End_k(M_Q)$ as $k$-algebras and as $(\Br_R(A),\Br_R(A))$-bimodules.

\end{enumerate}
\end{proposition}

\begin{proposition}[{\cite[Proposition 4.1]{Lin15}}]\label{lin154.1}
Let $A$ and $B$ be almost source algebras of blocks of finite group algebras having a common defect group $P$ and the same fusion system $\mathcal{F}$ on $P$. Let $V$ be an $\mathcal{F}$-stable indecomposable endopermutation $kP$-module with vertex $P$. Let $M$ be a direct summand the $(A,B)$-bimodule
$$U:=A\otimes_{kP}{\rm Ind}_{\Delta P}^{P\times P}(V)\otimes_{kP}B.$$
Consider $M$ as a $k\Delta P$-module via the homomorphism $\Delta P\to A\otimes_k B^{\rm op}$ sending $(u,u)\in \Delta P$ to $u1_A\otimes u^{-1}1_B$. Then for any non-trivial subgroup $Q$ of $P$, there is a canonical $({\rm Br}_Q(A), \Br_Q(B))$-bimodule $M_Q$ satisfying ${\rm Br}_{\Delta Q}({\rm End}_k(M))\cong {\rm End}_k(M_Q)$ as $k$-algebras and as $(\Br_Q(A)\otimes_k\Br_Q(A)^{\rm op},\Br_Q(B)\otimes_k\Br_Q(B)^{\rm op})$-bimodules.
\end{proposition}


\begin{theorem}[a slight refinement of {\cite[Theorem 4.2]{Lin15}} or {\cite[Proposition 9.9.7]{Linckelmann}}]\label{lin154.2}
Keep the notation of Proposition \ref{lin154.1}. Then $\End_{B^{\rm op}}(U)$ is a $\Delta Q$-subalgebra of $\End_k(U)$, the algebra homomorphism
$$\beta:{\rm Br}_{\Delta Q}(\End_{B^{\rm op}}(U))\to {\rm Br}_{\Delta Q}(\End_k(U))$$
induced by the inclusion $\End_{B^{\rm op}}(U)\to \End_k(U)$ is injective, and there is a commutative diagram of algebra and also $(\Br_Q(A),\Br_Q(A))$-bimodule homomorphisms
$$\xymatrix{\Br_{\Delta Q}(\End_k(U))\ar[r]^{\cong} & \End_k(U_Q)\\
\Br_{\Delta Q}(\End_{B^{\rm op}}(U))\ar[u]^{\beta}\ar[r]^{\cong}_{\gamma} & \End_{\Br_Q(B)^{\rm op}}(U_Q)\ar[u]
}$$
where the right vertical arrow is the obvious inclusion map. Moreover, the homomorphism $\gamma$ is an isomorphism.
\end{theorem}

\noindent{\it Proof.} All of the statements are proved in \cite{Lin15} (or \cite{Linckelmann}) except for that $\gamma$ is surjective. So we borrow the proof of \cite[Theorem 4.2]{Lin15} (or \cite[Proposition 9.9.7]{Linckelmann}) to here and continue to prove that $\gamma$ is surjective. Denote by $Y$ the set of elements in ${\rm Br}_{\Delta Q}(\End_{B^{\rm op}}(U))$ which commute with all elements in $\Br_Q(B)$, then $Y$ is a $(\Br_Q(A),\Br_Q(A))$-bimodule. We see that the image of $\beta$ is contained in $Y$. Since the upper horizontal map is a $(\Br_Q(A)\otimes_k \Br_Q(A)^{\rm op},\Br_Q(B)\otimes_k \Br_Q(B)^{\rm op})$-bimodule homomorphism, it restricts an isomorphism $Y\cong \End_{\Br_Q(B)^{\rm op}}(U_Q)$ of $(\Br_Q(A),\Br_Q(A))$-bimodules. To show that $\gamma$ is surjective, it suffices to show that ${\rm Im}(\beta)=Y$. In order to show that ${\rm Im}(\beta)=Y$, we first note that this does not make use of the left $A$-module structure of $U$ but only of the left $kQ$-module structure. By the end of the proof of \cite[Theorem 4.2]{Lin15} (or \cite[Proposition 9.9.7]{Linckelmann}), it suffices to consider the structural homomorphism
$$\Br_Q(B)\to {\rm End}_k(\Br_Q(B))$$
of the right $\Br_Q(B)$-module $\Br_Q(B)$ instead of $\beta$.
The set of elements in ${\rm End}_k(\Br_Q(B))$ which commute with all elements in $\Br_Q(B)$ is ${\rm End}_{\Br_Q(B)^{\rm op}}(\Br_Q(B))$. So it suffices to show that the map
$$\Br_Q(B)\to {\rm End}_{\Br_Q(B)^{\rm op}}(\Br_Q(B))$$
is surjective. But this map is obviously an isomorphism. $\hfill\square$

\subsection{Brauer-friendly modules and slashed modules}

A Brauer-friendly module is a direct sum of indecomposable modules with compatible fusion-stable endopermutation sources. We refer to \cite[Definition 8]{Biland} for its definition.

\begin{proposition}\label{brauerfriendly}
	Let $G$ and $H$ be finite groups, and let $b$ and $c$ be blocks of $kG$ and $kH$ respectively. Let $M$ be an indecomposable
	$(kGb, kHc)$-bimodule with an endopermutation source inducing a stable equivalence of Morita type between $kGb$ and $kHc$. Then $M$ is a Brauer-friendly $k(G\times H)(b\otimes c^\circ)$-module.
\end{proposition}

To prove this proposition, we first prove the following lemma.

\begin{lemma}\label{source triple}
	With the notation and assumptions from Theorem \ref{Puig1} and Theorem \ref{main}, $(\Delta P, e_P\otimes f_P^{\circ}, V)$ is a source triple of $M$ in the sense of \cite[Definition 2]{Biland}.
\end{lemma}

\noindent{\it Proof.} By \cite[Lemma 1 (ii)]{Biland}, it suffices to prove that the $k\Delta P$-module $e_PMf_P$ admits an indecomposable direct summand which is isomorphic to $V$. We first show that the $k\Delta P$-module $e_PMf_P$ admits an indecomposable direct summand with vertex $\Delta P$. Indeed, we have
\begin{align*}
	{\rm Br}_{\Delta P}(\End_k(e_PMf_P))&\cong \Br_{\Delta P}((e_P\otimes f_P^\circ)(M^*\otimes_k M)(e_P\otimes f_P^\circ))\\
	&\cong (e_P\otimes f_P^\circ)\Br_{\Delta P}(M^*\otimes_k M)(e_P\otimes f_P^\circ)
\end{align*}
(at least) as $k$-modules, where the second isomorphism is by \cite[Lemma 3.9]{Lin08}. By Proposition \ref{slash2} (i), the $k\Delta P$-module $iMj$ has a direct summand isomorphic to $V$ which has the vertex $\Delta P$. So we have
\begin{align*}
	0\neq{\rm Br}_{\Delta P}(\End_k(iMj))&\cong \Br_{\Delta P}((i\otimes j^\circ)(M^*\otimes_k M)(i\otimes j^\circ))\\
	&\cong {\rm br}_{\Delta P}^{k(G\times H)}(i\otimes j^\circ)\Br_{\Delta P}(M^*\otimes_k M){\rm br}_{\Delta P}^{k(G\times H)}(i\otimes j^\circ),
\end{align*}
(at least) as $k$-modules,  where the second isomorphism is again by \cite[Lemma 3.9]{Lin08}.
Since (by \cite[Proposition 4.9]{Lin08})
$$(e_P\otimes f_P^\circ){\rm br}_{\Delta P}^{k(G\times H)}(i\otimes j^\circ)={\rm br}_{\Delta P}^{k(G\times H)}(i\otimes j^\circ)={\rm br}_{\Delta P}^{k(G\times H)}(i\otimes j^\circ)(e_P\otimes f_P),$$
we deduce that ${\rm Br}_{\Delta P}(\End_k(e_PMf_P))\neq 0$. Then by \cite[Corollary 2.6.8]{Linckelmann}, the $k\Delta P$-module $e_PMf_P$ admits at least an indecomposable direct summand with vertex $\Delta P$.

By Theorem \ref{Puig1}, $M$ is isomorphic to a direct summand of $k(G\times H)(i\otimes j^\circ)\otimes_{k\Delta P} V$, hence $e_PMf_P$ is isomorphic to a direct summand of $(e_P\otimes f_P^\circ)k(G\times H)(i\otimes j^\circ)\otimes_{k\Delta P} V$. Using a similar argument as in the proof of \cite[Proposition 5.2 (i)]{Lin08}, we can easily prove that as a $(k\Delta P,k\Delta P)$-bimodule, every indecomposable direct summand of $(e_P\otimes f_P^\circ)k(G\times H)(i\otimes j^\circ)$ is isomorphic to $k\Delta P\otimes_{k\Delta S} {}_\varphi k\Delta P$ for some subgroup $S$ of $P$ and some morphism $\varphi:\Delta S\to \Delta P$ belonging to $\mathcal{F}\times\mathcal{F}$.
It follows that every indecomposable direct summand of $e_PMf_P$ is isomorphic to a $k\Delta P$-module of the form $k\Delta P\otimes_{k\Delta S}{}_\varphi V$. This implies that $e_PMf_P$ has at least a direct summand isomorphic to ${}_\varphi V$ for some $\varphi\in {\rm Aut}_{\F\times\F}(\Delta P)$. By Theorem \ref{Puig1} (b), we have ${}_\varphi V\cong V$, which completes the proof.
$\hfill\square$

\medskip\noindent{\it Proof of Proposition \ref{brauerfriendly}.} Since we are in the context of Theorem \ref{Puig1}, we can use the notation and assumptions from Theorem \ref{Puig1} and Theorem \ref{main}. By Lemma \ref{source triple}, $(\Delta P, e_P\otimes f_P^{\circ}, V)$ is a source triple of $M$ in the sense of \cite[Definition 2]{Biland}. Since $M$ is indecomposable, by \cite[Definition 8]{Biland}, it suffices to show that the source triple $(\Delta P, e_P\otimes f_P^{\circ}, V)$ is fusion-stable in the sense of \cite[Definition 5]{Biland}. By the paragraph below \cite[Definition 5]{Biland}, it suffices to show that ${\rm Res}_{\Delta Q}^{\Delta P}(V)\oplus {}_\psi V$ is an endopermutation $k\Delta Q$-module, for any $(b\otimes c^{\circ})$-subpair $(\Delta Q, e_Q\otimes f_Q^{\circ})$ contained in $(\Delta P, e_P\otimes f_P^{\circ})$ and any group homomorphism $\psi: \Delta Q\to \Delta P$ induced by conjugation by an element $(g,h)\in G\times H$ satisfying ${}^{(g,h)}(\Delta Q,e_Q\otimes f_Q^{\circ})\subseteq(\Delta P,e_P\otimes f_P^{\circ})$. Note that $\Delta P$ and $\Delta Q$ are objects of the fusion system $\F\times\mathcal{F}$ of the block $b\otimes c^{\circ}$ determined by the almost idempotent $i\otimes j^{\circ}$ (see \cite[Proposition 4.9]{Lin08}). Since $\F\times\F$ is a full subcategory of the Brauer category (in the sense of \cite{Biland}) of the block $b\otimes c^{\circ}$, $\psi$ is a morphism in $\mathcal{F}\times \mathcal{F}$.
So it suffices to show that ${\rm Res}_Q^P(V)\oplus {}_\varphi V$ is an endopermutation $kQ$-module for any $b$-subpair $(Q,e_Q)$ contained in $(P,e_P)$ and any morphism $\varphi: Q\to P$ in $\F$. This is ensured by Theorem \ref{Puig1} (b). $\hfill\square$

\medskip Biland showed that Dade's slash construction for endopermutation modules, also known as deflation-restriction, can be extended for applying to Brauer-friendly modules.

\begin{lemma}[see {\cite[Theorem 18]{Biland}} or {\cite[Lemma 3]{Bilandadv}}]\label{slash}
	Let $G$ be a finite group, $b$ a block of $kG$, $M$ a Brauer-friendly $kGb$-module, $(P,e_P)$ a $b$-subpair, and $H$ a subgroup of $G$ such that $PC_G(P)\leq H\leq N_G(P,e_P)$. Then there exists a Brauer-friendly $kHe_P$-module $Sl_{(P,e_P)}^H(M)$ and an isomorphism of $C_G(P)$-interior $H$-algebras
	$$\theta_{(P,e_P)}^H:~{\rm Br}_P({\rm End}_k(e_PM))\cong {\rm End}_k(Sl_{(P,e_P)}^H(M)).$$
\end{lemma}

Following Biland, the pair $(Sl_{(P,e_P)}^H(M),\theta_{(P,e_P)}^H)$ or just the $kHe_P$-module $Sl_{(P,e_P)}^H(M)$ is called a $(P,e_P)$-{\it slashed module} attached to $M$ over the group $H$. If $M$ is a $p$-permutation module, then there is a canonical choice of $(P,e_P)$-slashed module attached to $M$ over $H$: the Brauer quotient ${\rm Br}_P(e_PM)$ (regarded as a $kHe_P$-module via restriction), together with the natural isomorphism ${\rm Br}_P({\rm End}_k(e_PM))\cong {\rm End}_k({\rm Br}_P(e_PM))$ (see \cite[Lemma 3.3]{Broue85}).

\section{Proof of Theorem \ref{main2}}

\begin{proposition}\label{Brauerconstruction}
Let $A$ be an almost source algebra of a block a finite group algebra over $k$ with defect group $P$, and $Q$ a subgroup of $P$. If $M$ is a projective $A$-module, then $M$ is also projective as a $kQ$-module. In particular, if $Q\neq 1$, ${\rm Br}_Q(M)=0$.
\end{proposition}
\noindent{\it Proof.} It suffices to show that $A$ is projective as a $kQ$-module. As a $kQ$-module, $A$ is isomorphic to a direct summand of ${\rm Res}_{Q}^{G}(kG)= \oplus_{g\in [Q\setminus G]} kQg$, hence projective. $\hfill\square$

\begin{proposition}\label{dualsourcealgebra}
Keep the notation and assumptions in Proposition \ref{lin154.1}. Then $M_Q^*$ is a $({\rm Br}_Q(B), \Br_Q(A))$-bimodule satisfying ${\rm End}_k(M_Q^*)\cong {\rm Br}_{\Delta Q}({\rm End}_k(M^*))$ as $k$-algebras and as $({\rm Br}_Q(B)\otimes_k{\rm Br}_Q(B)^{\rm op}, \Br_Q(A)\otimes_k {\rm Br}_Q(A)^{\rm op})$-bimodules.
\end{proposition}

\noindent{\it Proof.}  We have $${\rm Br}_{\Delta Q}({\rm End}_k(M^*))\cong{\rm Br}_{\Delta Q}({\rm End}_k(M)^{\rm op})\cong \Br_{\Delta Q}(\End_k(M))^{\rm op}\cong{\rm End}_k(M_Q)^{\rm op}\cong{\rm End}_k(M_Q^*)$$
as $k$-algebras and as $({\rm Br}_Q(B)\otimes_k{\rm Br}_Q(B)^{\rm op}, \Br_Q(A)\otimes_k {\rm Br}_Q(A)^{\rm op})$-bimodules,
where the first and fourth isomorphisms are by \cite[Proposition 2.9.4]{Linckelmann} and \cite[\S4.1.2]{CR}, and the third isomorphism is by the definition of $M_Q$.
  $\hfill\square$

\medskip\noindent{\it Proof of Theorem \ref{main2}.}
By Theorem \ref{lin154.2}, we have an isomorphism
$$\Br_{\Delta Q}(\End_{B^{\rm op}}(M))\to {\rm End}_{\Br_Q(B)^{\rm op}}(M_Q)$$
of algebras and $(\Br_Q(A),\Br_Q(A))$-bimodules.
By \cite[Corollary 2.12.4]{Linckelmann}, $\End_{B^{\rm op}}(M)\cong M\otimes_B M^*$ as $(A,A)$-bimodules and ${\rm End}_{\Br_Q(B)^{\rm op}}(M_Q)\cong M_Q\otimes_{\Br_Q(B)} M_Q^*$ as $(\Br_Q(A),\Br_Q(A))$-bimodules. Since $M$ and $M^*$ induce a stable equivalence of Morita type between $A$ and $B$, $M\otimes_B M^*\cong A\oplus W$ for some projective $A\otimes_k A^{\rm op}$-module $W$. Since ${\rm Br}_{\Delta Q}(W)=0$ (see Proposition \ref{Brauerconstruction}),
$$\Br_{\Delta Q}(\End_{B^{\rm op}}(M))\cong\Br_{\Delta Q}(A\oplus W)\cong \Br_Q(A)$$
as $(\Br_Q(A),\Br_Q(A))$-bimodules.
So we have an isomorphism
$${\rm Br}_Q(A)\cong  M_Q\otimes_{\Br_Q(B)} M_Q^*$$
of $(\Br_Q(A),\Br_Q(A))$-bimodules. Applying the same argument to the bimodules $M^*$ and $M_Q^*$ with the roles of $A$ and $B$ exchanged, we can prove that $${\rm Br}_Q(B)\cong  M_Q^*\otimes_{\Br_Q(A)} M_Q$$
as $(\Br_Q(B),\Br_Q(B))$-bimodules.  $\hfill\square$

\section{Proofs of Theorem \ref{main} and Corollary \ref{coro1}}

The following proposition is an analogue of Proposition \ref{dualsourcealgebra} at the block algebra level.
\begin{proposition}\label{dual}
Let $M$ be a Brauer-friendly $kGb$-module, $(P,e_P)$ a $b$-subpair, $H$ a subgroup of $G$ satisfying $PC_G(P)\leq H\leq N_G(P,e_P)$, and $M_P$ a $(P,e_P)$-slashed module attached to $M$ over $H$. Then $M^*$ is a Brauer-friendly $kGb^\circ$-module, and $M_P^*$ is a $(P,e_P^\circ)$-slashed module attached to $M^*$ over $H$.
\end{proposition}

\noindent{\it Proof.} Let $M=\bigoplus_{i=1}^sM_i$ be a decomposition of $M$ into indecomposable $kGb$-module. Then each $M_i^*$ is an indecomposable $kGb^\circ$-module, and $M^*=\bigoplus_{i=1}^sM_i^*$. Let $(Q_i,e_i,V_i)$ be a source triple of $M_i$. (We refer to \cite[Definition 2]{Biland} for the definition of source triples.) By definition, it is easy to see that $(Q_i,e_i^\circ,V_i^*)$ is a source triple of $M_i$. For any $i,j\in \{1,\cdots,s\}$, since $(Q_i,e_i,V_i)$ and $(Q_j,e_j,V_j)$ are compatible fusion-stable endopermutation source triples, it is easy to check that $(Q_i,e_i^\circ,V_i^*)$ and $(Q_j,e_j^\circ,V_j^*)$ are also compatible fusion-stable endopermutation source triples. Hence $M^*$ is a Brauer-friendly $kGb^\circ$-module.

Noting that $e_P^\circ M^*$ is the dual of $e_PM$, we have isomorphisms of $C_G(P)$-interior $H$-algebras
$${\rm Br}_{P}({\rm End}_k(e_P^\circ M^*))\cong{\rm Br}_{P}({\rm End}_k(e_PM)^{\rm op})\cong{\rm Br}_{P}({\rm End}_k(e_PM))^{\rm op} \cong{\rm End}_k(M_P)^{\rm op}\cong{\rm End}_k(M_P^*).$$
(See e.g. \cite[Proposition 2.9.4]{Linckelmann} and \cite[\S4.1.2]{CR} for the first and fourth isomorphisms.) So by definition, $M_P^*$ is a $(P,e_P^\circ)$-slashed module attached to $M^*$ over $H$.
$\hfill\square$

\medskip\noindent{\it Proof of Theorem \ref{main}.} (a). Let $A:=ikGi$, $B:=jkHj$ and $U:=iMj$. By the definition of a $(\Delta Q, e_Q\otimes f_Q^\circ)$-slashed module attached to $M$ over $N_{G\times H}(\Delta Q,e_Q\otimes f_Q^\circ)$, after restriction, $M_Q$ is a $(kC_G(Q)e_Q,kC_H(Q)f_Q)$-bimodule such that $\End_k(M_Q)\cong \Br_{\Delta Q}({\rm End}_k(e_QMf_Q))$ as algebras and as $(kC_G(Q)e_Q\otimes_k (kC_G(Q)e_Q)^{\rm op},kC_H(Q)f_Q\otimes_k (kC_H(Q)f_Q)^{\rm op})$-bimodules. We first consider the case that $Q$ is fully $\F$-centralised. By Proposition \ref{almost}, ${\rm br}_Q^{kG}(i)$ and ${\rm br}_Q^{kH}(j)$ are almost source idempotents of $kC_G(Q)e_Q$ and $kC_H(Q)f_Q$, respectively. Let $U_Q:={\rm br}_Q^{kG}(i)M_Q{\rm br}_Q^{kH}(j)$. By the standard Morita equivalences between block algebras and almost source algebras, $U_Q$ is a $({\rm Br}_Q(A), \Br_Q(B))$-bimodule such that ${\rm End}_k(U_Q)\cong {\rm Br}_{\Delta Q}({\rm End}_k(U))$ as algebras and as $(\Br_Q(A)\otimes_k\Br_Q(A)^{\rm op},\Br_Q(B)\otimes_k\Br_Q(B)^{\rm op})$-bimodules. By Theorem \ref{main2}, $U_Q$ induces a Morita equivalence between ${\rm Br}_Q(A)$ and $\Br_Q(B)$. Hence $M_Q$ induces a Morita equivalence between $kC_G(Q)e_Q$ and $kC_H(Q)f_Q$.

Now we consider general $Q$. We can always find a subgroup $R$ of $P$ such that $Q\cong R$ in $\F$ and such that $R$ is fully $\F$-centralised. Assume that the isomorphism $Q\cong R$ is induced by an element $g\in G$ (resp. an element $h\in H$). Namely that the isomorphism $Q\cong R$ sends $u\in Q$ to $gug^{-1}=huh^{-1}$. By \cite[Lemma 22 (ii)]{Biland}, the $kN_{G\times H}(\Delta R, e_R\otimes f_R^\circ)(e_R\otimes f_R^\circ)$-module ${}^{(g,h)}M_Q$ is a $(\Delta R, e_R\otimes f_R^\circ)$-slashed module attached to $M$. By \cite[Lemma 17 (ii) or Theorem 18 (ii)]{Biland}, ${}^{(g,h)}M_Q\cong M_R$ as $(kC_G(R)e_R, kC_H(R)f_R)$-bimodules. So ${}^{(g,h)}M_Q$ induces a Morita equivalence between $kC_G(R)e_R$ and $kC_H(R)f_R$. This is equivalent to that $M_Q$ induces a Morita equivalence between $kC_G(Q)e_Q$ and $kC_H(Q)f_Q$.

Since ${\rm Res}_{C_G(Q)\times C_H(Q)}(M_Q)$ is a $(\Delta Q,e_Q\otimes f_Q^\circ)$-slashed module attached to $M$ over $C_G(Q)\times C_H(Q)$, it is Brauer-friendly (see Lemma \ref{slash} or \cite[Theorem 21]{Biland}), and hence has an endopermutation module as a source. Alternatively, since $M_Q$ is a Brauer-friendly $kN_{G\times H}(\Delta Q, e_Q\otimes f_Q^\circ)(e_Q\otimes f_Q^\circ)$-module (see Lemma \ref{slash} or \cite[Theorem 21]{Biland}), using the Mackey formula, it is easy to see that as a $(kC_G(Q)e_Q, kC_H(Q)f_Q)$-bimodule, $M_Q$ has an endopermutation module as a source.

\noindent(b). Since $(a)$ is proved, by \cite[Lemma 10.2.8]{Rou}, ${\rm Ind}_{T_Q}^{G_Q\times H_Q}(M_Q)$ induces a Morita equivalence between $kG_Qe_Q$ and $kH_Qf_Q$. Using the Mackey formula, we see that ${\rm Res}_{T_Q}(M_Q)$ has an endopermutation module as a source. By \cite[Proposition 5.1.7]{Linckelmann}, a source of ${\rm Res}_{T_Q}(M_Q)$ remains a source of ${\rm Ind}_{T_Q}^{G_Q\times H_Q}(M_Q)$.  $\hfill\square$

\medskip\noindent{\it Proof of Corollary \ref{coro1}.} Since $(\Delta P,e_P\otimes e_P^\circ, V)$ is a source triple of the $(kGb,kHc)$-bimodule $M$, by \cite[Theorem 23]{Biland}, $M_P$ is an indecomposable projective $kN_{G\times H}(\Delta P, e_P\otimes f_P^\circ)/\Delta P$-module, hence a $p$-permutation $kN_{G\times H}(\Delta P, e_P\otimes f_P^\circ)$-module. It follows that ${\rm Res}_{C_G(P)\times C_H(P)}^{N_{G\times H}(\Delta P, e_P\otimes f_P^\circ)}(M_P)$ and ${\rm Ind}_{T_P}^{G_P\times H_P}(M_P)$ are $p$-permutation modules.  This proves (a) and (b). Taking $G_P=N_G(P,e_P)$ and $H_P=N_H(P,f_P)$, the statement (c) follows by (b) and \cite[Theorem 6.2.6 (iii)]{Linckelmann}.  $\hfill\square$

\section{Proof of Corollary \ref{functorialequ}}

Let $\mathbb{F}$ be an algebraically closed field of characteristic 0. We refer to \cite[Definition 6.1]{BY} for the definition of a diagonal $p$-permutation functor over $\mathbb{F}$. Diagonal $p$-permutation functors over $\mathbb{F}$ form an abelian category which is denoted by $\F_{\mathbb{F}pp_k}^\Delta$. According to \cite[Corollary 6.15]{BY}, the category $\F_{\mathbb{F}pp_k}^\Delta$ is semisimple, and moreover, the simple diagonal $p$-permutation functors, up to isomorphism, are parametrized by the isomorphism classes of triples
$(L,u,V)$ where $(L,u)$ is a $D^\Delta$-pair, and $V$ is a simple $\mathbb{F}{\rm Out}(L,u)$-module. We refer to \cite[Definition 3.2]{BY} for the definition of a $D^\Delta$-pair and refer to \cite[Notation 6.8]{BY} for the definitions of the groups ${\rm Aut}(L,u)$ and ${\rm Out}(L,u)$.

Let $G$ be a finite group and $b$ a block of $kG$. The pair $(G,b)$ defines an object $\mathbb{F}T_{G,b}^\Delta$ of the category $\F_{\mathbb{F}pp_k}^\Delta$, which is called a block diagonal $p$-permutation functor (see the beginning of \cite[\S8]{BY}). Let $H$ be another finite group and $c$ a block of $kH$, then the pairs $(G,b)$ and $(H,c)$ are defined to be {\it functorially equivalent} over $\mathbb{F}$ if the objects $\mathbb{F}T_{G,b}^{\Delta}$ and $\mathbb{F}T_{H,c}^\Delta$ are isomorphic in $\F_{\mathbb{F}pp_k}^\Delta$ (see \cite[Definition 10.1]{BY}).

\begin{remark}\label{remarkSLuV}
{\rm \begin{enumerate} [\rm (i)]
		
\item		Since the category $\F_{\mathbb{F}pp_k}^\Delta$ is semisimple, the functor $\mathbb{F}T_{G,b}^\Delta$ is a direct sum of simple diagonal $p$-permutation functors $S_{L,u,V}$. Hence the pairs $(G,b)$ and $(H,c)$ are functorially equivalent over $\mathbb{F}$ if and only if for any triple $(L,u,V)$, the multiplicities of the simple diagonal $p$-permutation functor
$S_{L,u,V}$ in $\mathbb{F}T_{G,b}^\Delta$ and $\mathbb{F}T_{H,c}^\Delta$ are the same.

\item In \cite{BY2}, Bouc and Y{\i}lmaz defined the notion of stable functorial equivalences over $\mathbb{F}$, and they (\cite[Theorem 1.2 (i)]{BY2}) proved that the pairs $(G,b)$ and $(H,c)$ are stably functorially equivalent over $\mathbb{F}$ if and only if for any triple $(L,u,V)$ with $L\neq 1$, the multiplicities of
$S_{L,u,V}$ in $\mathbb{F}T_{G,b}^\Delta$ and $\mathbb{F}T_{H,c}^\Delta$ are the same. 
\end{enumerate}
}
\end{remark}

Hence we now review the formula (given in \cite[Theorem 8.22 (b)]{BY}) for the multiplicity of $S_{L,u,V}$ in $\mathbb{F}T_{G,b}^\Delta$. Let $(P,e_P)$ be a maximal $b$-subpair. For any subgroup $Q\leq P$, let $e_Q$ be the unique block of $kC_G(Q)$ with $(Q, e_Q)$ contained in $(P, e_P)$. Let $\F$ be the fusion system of $b$ with respect to $(P,e_P)$ (see \cite[Definition 8.5.1]{Linckelmann}) and let
$[\F]$ be a set of isomorphism classes of objects in $\F$. For $Q\in \F$, denote by $\mathcal{P}_{(Q,e_Q)}^G(L,u)$ the set of group isomorphisms $\pi:L\to Q$ with $\pi i_u \pi^{-1}\in \Aut_{\mathcal{F}}(Q)$, where $i_u$ denotes the conjugation by $u$. The set $\mathcal{P}_{(Q,e_Q)}^G(L,u)$ is an $(N_G(Q,e_Q), {\rm Aut}(L,u))$-biset via
\begin{align*}
	g\cdot \pi\cdot\alpha = i_g\pi\alpha
\end{align*}
for $g\in N_G(Q,e_Q)$, $\pi\in \mathcal{P}_{(Q,e_Q)}^G(L,u)$ and $\alpha\in{\rm Aut}(L,u)$. Denote by $[\mathcal{P}_{(Q,e_Q)}^G(L,u)]$ a set of representatives of $N_G(Q,e_Q)\times \Aut(L,u)$-orbits of $\mathcal{P}_{(Q,e_Q)}^G(L,u)$.  For $\pi \in [\mathcal{P}_{(Q,e_Q)}^G(L,u)]$, we denote by ${\rm Aut}(L,u)_{\overline{(Q,e_Q,\pi)}}$ the stabiliser in ${\rm Aut}(L,u)$ of the $N_G(Q,e_Q)$-orbit of $\pi$. Then we have
\begin{align*}
	{\rm Aut}(L,u)_{\overline{(Q,e_Q,\pi)}}=\{\alpha\in {\rm Aut}(L,u)\mid \pi\alpha \pi^{-1}\in {\rm Aut}_{\F}(Q)\}\,.
\end{align*}

Denote by $\Proj(kGb)$ the set of isomorphism classes of indecomposable projective $kGb$-modules. Let $\varphi$ be an automorphism of the group $G$. We use abusively the same notation $\varphi$ to denote the restriction of $\varphi$ to any subgroup of $G$. We also use the same notation $\varphi$ to denote the obvious $k$-algebra automorphism of $kG$ induced by $\varphi$. For a $kG$-module $M$, denote by ${}_\varphi M$ the $kG$-module which is equal $M$ as a $k$-module endowed with the structure homomorphism $G\xrightarrow{\varphi}G\to {\rm End}_k(M)$. If $\varphi(b)=b$, we also use abusively the same notation $\varphi$ to denote bijection $\Proj(kGb)\to \Proj(kGb)$ sending $[S]\in \mathcal{P}$ to $[{}_\varphi S]$. For a $G$-algebra $A$, we denote by ${}_\varphi A$ the $G$-algebra which is equal to $A$ as a $k$-algebra endowed with the structure homomorphism $G\xrightarrow{\varphi}G\to {\rm Aut}(A)$.

\begin{lemma}\label{lemma:automorphism acts on modules}
Let $\varphi$ and $\psi$ be elements in $\Aut(G)$ satisfying $\varphi(b)=b$ and $\psi(b)=b$. Let $M$ be a $kGb$-module. If ${}_\varphi M\cong M$, then ${}_{\varphi\circ\psi}M\cong {}_\psi M$. 
\end{lemma}

\noindent{\it Proof.}  Let $f$ be an isomorphism from ${}_\varphi M$ to $M$. Hence for any $a\in kGb$ and $m\in M$, we have $f(\varphi(a)m)=af(m)$. We claim that $f$ is also an isomorphism from ${}_{\varphi\circ\psi}M$ to ${}_\psi M$. Indeed, for any $a\in kGb$ and $m\in M$, we have
$$f((\varphi\circ\psi)(a)m)=f(\varphi(\psi(a))m)=\psi(a)f(m),$$
and this implies the claim. $\hfill\square$

\medskip Let $\pi:L\to Q$ be an element in $\mathcal{P}_{(Q,e_Q)}^G(L,u)$. That is, $\pi i_u\pi^{-1}=i_g\in \Aut_\F(Q)$ for some $g\in N_G(Q,e_Q)$. We say that a $kC_G(Q)e_Q$-module $M$ is {\it $u$-invariant with respect to $\pi$} if ${}_{i_g}M\cong M$. Let $$\Proj(kC_G(Q)e_P, u,\pi):=\{[S]\in \Proj(kC_G(Q)e_Q)~|~\mbox{$S$ is $u$-invariant with respect to $\pi$}\}.$$

\begin{remark}
{\rm In \cite{BY}, $\Proj(kC_G(Q)e_Q, u,\pi)$ is written as $\Proj(kC_G(Q)e_Q,u)$ without the reference of $\pi$.}
\end{remark}	

For any $\alpha\in \Aut(L,u)_{\overline{(Q,e_Q,\pi)}}$, we have $\pi\alpha\pi^{-1}=i_{g'}\in \Aut_\F(Q)$ for some $g'\in N_G(Q,e_Q)$. Let $[S]\in \Proj(kC_G(Q)e_Q,u,\pi)$, then by Lemma \ref{lemma:automorphism acts on modules}, we have ${}_{i_g\circ i_{g'}}S\cong{}_{i_{g'}}S$, hence $[{}_{i_{g'}}S]\in \Proj(kC_G(Q)e_Q,u,\pi)$. One checks that the formula $[S]\cdot \alpha:=[{}_{i_{g'}^{-1}}S]$ defines a right action of $\Aut(L,u)_{\overline{(Q,e_Q,\pi)}}$ on $\Proj(kC_G(Q)e_Q,u,\pi)$. This right action defines a right permutation $\mathbb{F}\Aut(L,u)_{\overline{(Q,e_Q,\pi)}}$-module, which we denote it by $\mathbb{F}\Proj(kC_G(Q)e_Q,u,\pi)$.

\begin{theorem}[{\cite[Theorem 8.22 (b)]{BY}}]\label{theorem:BY,multiplicity}
The multiplicity of a simple diagonal $p$-permutation functor $S_{L,u,V}$ in the functor $\mathbb{F}T_{G,b}^\Delta$ is equal to the $\mathbb{F}$-dimension of the $\mathbb{F}$-vector space
\begin{align*}
	\bigoplus_{Q \in [\F]} \bigoplus_{\pi \in [\mathcal{P}_{(Q,e_Q)}^G(L,u)]} \mathbb{F}{\rm Proj}(kC_G(Q)e_Q,u,\pi)\otimes_{\FF{\rm Aut}(L,u)_{\overline{(Q,e_Q,\pi)}}} V\,.
\end{align*}
\end{theorem}

\begin{lemma}[{a slight modification of \cite[Lemma 3.2]{Yilmaz}}]\label{lem:criterion}
Let $G$ and $H$ be finite groups, $b$ a block of $kG$ and $c$ a block of $kH$.  Let $(D,e_D)$ and $(E,f_E)$ be maximal $b$- and $c$- subpairs, respectively. For any subgroup $P\le D$, let $e_P$ denote the unique block of $kC_G(P)$ with $(P,e_P)\subseteq (D,e_D)$ and for any subgroup $Q\le E$, let $f_Q$ denote the unique block of $kC_H(Q)$ with $(Q,f_Q)\subseteq (E,f_E)$.  Let $\F$ and $\G$ denote the fusion systems of $b$ and $c$ with respect to the maximal subpairs $(D,e_D)$ and $(E,f_E)$, respectively. Suppose that there is a group isomorphism $\phi: E\to D$ such that $\F={}^\phi\G$.  Suppose further that for every subgroup $Q\leq E$ (resp. $1< Q\leq E$), we have a bijection
\begin{align*}
	\Psi_Q: \Proj(kC_H(Q)f_Q)\to \Proj(kC_G(\phi(Q))e_{\phi(Q)})
\end{align*}
such that for any $(g,h)\in N_G(\phi(Q), e_{\phi(Q)})\times N_H(Q,f_Q)$ the maps $i_g\Psi_Q i_h$ and $\Psi_Q$ are the same. Then the pairs $(G,b)$ and $(H,c)$ are functorially (resp. stably functorially) equivalent over $\FF$.
	\end{lemma}
	
\noindent{\it Proof (a slight modification of the proof of \cite[Lemma 3.2]{Yilmaz}).} Using the formula in Theorem \ref{theorem:BY,multiplicity}, we will show that the multiplicity of any simple functors $S_{L,u,V}$ (resp. $S_{L,u,V}$ with $L\neq 1$) in $\FFTD_{G,b}$ and in $\FFTD_{H,c}$ are the same. Let $Q\in \G$ with $L\cong Q$ and set $P:=\phi(Q)\in \F$. First, note that the map
\begin{align*}
	\P_{(Q,f_Q)}^H(L,u)\to \P_{(P,e_P)}^G(L,u), \quad \pi\mapsto \phi\pi
\end{align*}
is a bijection. Indeed, if $\pi\in \P_{(Q,f_Q)}^H(L,u)$, then $\pi i_u\pi^{-1} =i_h\in\Aut_{\G}(Q)$ for some $h\in N_H(Q,f_Q)$. Since $\phi$ induces an isomorphism on fusion systems, it follows that $\phi \pi i_u\pi^{-1}\phi^{-1}=\phi i_h\phi^{-1}\in\Aut_{\F}(P)$ and hence $\phi\pi\in \P_{(P,e_P)}^G(L,u)$. 

Now assume that two maps $\pi,\rho \in \P_{(Q,f_Q)}^H(L,u)$ are in the same $N_H(Q,f_Q)\times \Aut(L,u)$-orbit. Let $h\in N_H(Q,f_Q)$ and $\alpha\in\Aut(L,u)$ with
\begin{align*}
	\rho=i_h \pi\alpha\,.
\end{align*}
Let $g\in N_G(P,e_P)$ with $\phi i_h\phi^{-1}=i_g$. Then we have
\begin{align*}
	\phi\rho=\phi i_h\pi\alpha = \phi i_h\phi^{-1}\phi \pi\alpha =i_g \phi\pi \alpha\,.
\end{align*}
This shows that $\phi\rho$ and $\phi\pi$ are in the same $N_G(P,e_P)\times \Aut(L,u)$-orbit. Hence we also have a bijection between the sets $[\P_{(Q,f_Q)}^H(L,u)]$ and $[\P_{(P,e_P)}^G(L,u)]$ of double coset representatives.

Let $\pi\in [\P_{(Q,f_Q)}^H(L,u)]$. We have
\begin{align*}
	\Aut(L,u)_{\overline{(Q,f_Q,\pi)}}&=\{\alpha\in\Aut(L,u)\mid \exists h\in N_H(Q,f_Q),\, i_h=\pi\alpha\pi^{-1}\}\\&
	=\{\alpha\in\Aut(L,u)\mid \exists h\in N_H(Q,f_Q),\, \phi i_h\phi^{-1}=\phi\pi\alpha\pi^{-1}\phi^{-1}\}\\&
	=\{\alpha\in\Aut(L,u)\mid \exists g\in N_G(P,e_P),\, i_g=\phi\pi\alpha(\phi\pi)^{-1}\}\\&
	=\Aut(L,u)_{\overline{(P,e_P,\phi\pi)}}\,.
\end{align*}

Next we show that for $\pi\in [\P_{(Q,f_Q)}^H(L,u)]$ and hence for $\phi\pi\in [\P_{(P,e_P)}^G(L,u)]$, we have an isomorphism 
\begin{align*}
	\Proj(kC_H(Q)f_Q,u,\pi)\cong \Proj(kC_G(P)e_P,u,\phi\pi)
\end{align*}
of right $\Aut(L,u)_{\overline{(Q,f_Q,\pi)}}$-sets. Let $[S]\in \Proj(kC_H(Q)f_Q,u,\pi)$. Let also $h\in N_H(Q,f_Q)$ with $\pi i_u\pi^{-1}=i_h$. Then $[S]$ is $u$-invariant with respect to $\pi$ means that ${}_{i_h}S\cong S$, or equivalently, $i_h([S])=[S]$. Let $g\in N_G(P,e_P)$ with $\phi i_h\phi^{-1}=i_g$. Since $\Psi_Q=i_g\Psi_Q i_h$, it follows that 
\begin{align*}
	{i_g}(\Psi_Q([S]))=i_g\Psi_Q i_h([S])=\Psi_Q([S])
\end{align*}
and so $\Psi_Q([S])\in \Proj(kC_G(P)e_P,u,\phi\pi)$. This shows that the map $\Psi_Q$ restricts to a bijection
\begin{align*}
	\Psi_Q: \Proj(kC_H(Q)f_Q,u,\pi)\to \Proj(kC_G(P)e_P,u,\phi\pi).
\end{align*}
We show that this bijection is an isomorphism of right $\Aut(L,u)_{\overline{(Q,f_Q,\pi)}}$-sets. So let $\alpha\in \Aut(L,u)_{\overline{(Q,f_Q,\pi)}}$. Let also $h\in N_H(Q,f_Q)$ with $i_h =\pi\alpha\pi^{-1}$. Set again $\phi i_h\phi^{-1}=i_g$ for $g\in N_G(P,e_P)$. Let $[S]\in  \Proj(kC_H(Q)f_Q,u,\pi)$. Then,
\begin{align*}
	\Psi_Q([S]\cdot\alpha)=\Psi_Q([{}_{i_h}S])=\Psi_Q(i_h([S]))=i_g(\Psi_Q([S]))=\Psi_Q([S])\cdot \alpha\,.
\end{align*}

All these imply, by the formula in Theorem \ref{theorem:BY,multiplicity}, that the multiplicity of the simple diagonal $p$-permutation functor $S_{L,u,V}$ in $\FFTD_{G,b}$ and in $\FFTD_{H,c}$ are the same. Therefore, by Remark \ref{remarkSLuV}, the pairs $(G,b)$ and $(H,c)$ are functorially (resp. stably functorially) equivalent over $\FF$.  $\hfill\square$


\begin{proposition}\label{5.1}
Keep the notation and assumptions of Theorem \ref{main}. For every non-trivial subgroup $Q$ of $P$, we have a bijection
$$\Psi_Q:\Proj(kC_H(Q)f_Q)\to \Proj(kC_G(Q)e_Q)$$
such that for any $(g,h)\in N_G(Q,e_Q)\times N_H(Q,f_Q)$, the maps $i_g\Psi_Qi_h$ and $\Psi_Q$ are the same.
\end{proposition}

\noindent{\it Proof.} Since (by Theorem \ref{main} (a)) ${\rm Res}_{C_G(Q)\times C_H(Q)}(M_Q)$ induces a Morita equivalence between $kC_G(Q)e_Q$ and $kC_H(Q)f_Q$, the functor $M_Q\otimes_{kC_H(Q)f_Q}-$ induces a bijection
$$\Psi_Q:\Proj(kC_H(Q)f_Q)\to \Proj(kC_G(Q)e_Q).$$

The functor (say $F_g$) sending a $kC_G(Q)e_Q$-module $V$ to the $kC_G(Q)e_Q$-module ${}_{i_g}V$ induces a Morita self-equivalence of $kC_G(Q)e_Q$. This Morita equivalence is induced by the $(kC_G(Q)e_Q,kC_G(Q)e_Q)$-bimodule ${}_{i_g}(kC_G(Q)e_Q)$, which is equal to $kC_G(Q)e_Q$ as a right $kC_G(Q)e_Q$-module endowed as before the left action of $a\in kC_G(Q)e_Q$ by multiplication with $i_g(a)$. Similarly, the functor (say $F_h$) sending a $kC_H(Q)f_Q$-module $V$ to the $kC_H(Q)f_Q$-module ${}_{i_h}V$ induces a Morita self-equivalence of $kC_H(Q)f_Q$, which is induced by the $(kC_H(Q)f_Q,kC_H(Q)f_Q)$-bimodule ${}_{i_h}(kC_H(Q)f_Q)$. Hence the functor
$$F_g\circ(M_Q\otimes_{kC_H(Q)f_Q}-)\circ F_h$$
induces a Morita equivalence between $kC_G(Q)e_Q$ and $kC_H(Q)f_Q$, which is induced by the $(kC_G(Q)e_Q,kC_H(Q)f_Q)$-bimodule
$${}_{i_g}(kC_G(Q)e_Q)\otimes_{kC_G(Q)e_Q}M_Q\otimes_{kC_H(Q)f_Q}{}_{i_h}(kC_H(Q)f_Q)\cong {}_{i_g}(M_Q){}_{i_{h^{-1}}}.$$
This Morita equivalence in turn induces the bijection $i_g\Psi_Qi_h$.

Since $M$ is Brauer-friendly (see Proposition \ref{brauerfriendly}), by \cite[Lemma 7 (i)]{Biland}, it is easy to see that the $k(G\times H)(b\otimes c^{\circ})$-module ${}_{i_{g^{-1}}}M{}_{i_{h}}$ is a Brauer-friendly module.  Let $M_{g^{-1},h,Q}$ be a $(\Delta Q,e_Q\otimes f_Q^\circ)$-slashed module attached to ${}_{i_{g^{-1}}}M{}_{i_{h}}$ over $\Delta Q(C_G(Q)\times C_H(Q))$, we claim that $M_{g^{-1},h,Q}\cong {}_{i_g}(M_Q){}_{i_{h^{-1}}}$ as $(kC_G(Q)e_Q,kC_H(Q)f_Q)$-bimodules.
Note that ${}^{(g^{-1},h)}(\Delta Q)$ is the image of $\Delta Q$ under the automorphism $i_{(g^{-1},h)}$ of $G\times H$.
We have
$${\rm End}_k(M_{g^{-1},h,Q})\cong{\rm Br}_{\Delta Q}(\End_k(e_Q({}_{i_{g^{-1}}}M{}_{i_{h}})f_Q))={\rm Br}_{\Delta Q}(\End_k({}_{i_{(g^{-1},h)}}e_QMf_Q))$$
$$={\rm Br}_{\Delta Q}({}_{i_{(g^{-1},h)}}\End_k(e_QMf_Q))={\rm Br}_{{}^{(g^{-1},h)}(\Delta Q)}(\End_k(e_QMf_Q))\cong{}_{i_{(g,h^{-1})}}\Br_{\Delta Q}(\End_k(e_QMf_Q))$$
$$\cong {}_{i_{(g,h^{-1})}}\End_k(M_Q)=\End_k({}_{i_{(g,h^{-1})}}M_Q)=\End_k({}_{i_g}(M_Q){}_{i_{h^{-1}}})$$
as $(C_G(Q)\times C_H(Q))$-interior $\Delta Q(C_G(Q)\times C_H(Q))$-algebras. This means that ${}_{i_g}(M_Q){}_{i_{h^{-1}}}$ is another $(\Delta Q,e_Q\otimes f_Q^\circ)$-slashed module attached to ${}_{i_{g^{-1}}}M{}_{i_{h}}$ over $\Delta Q(C_G(Q)\times C_H(Q))$. Then by \cite[Theorem 18 (ii)]{Biland} or \cite[Lemma 3 (ii)]{Bilandadv} (actually, at this point, it is no need to use them), we have $M_{g^{-1},h,Q}\cong {}_{i_g}(M_Q){}_{i_{h^{-1}}}$ as $(kC_G(Q)e_Q,kC_H(Q)f_Q)$-bimodules, as claimed.

On the other hand, the $k(G\times H)(b\otimes c^{\circ})$-module ${}_{i_{g^{-1}}}M{}_{i_{h}}$ is isomorphic to $M$ via the map sending $m\in {}_{i_{g^{-1}}}M{}_{i_{h}}$ to $gmh\in M$. Using \cite[Theorem 18 (ii)]{Biland} again, we have ${}_{i_g}(M_Q){}_{i_{h^{-1}}}\cong M_Q$ as $(kC_G(Q)e_Q,kC_H(Q)f_Q)$-bimodules. So the maps $i_g\Psi_Qi_h$ and $\Psi_Q$ are the same. $\hfill\square$

\begin{remark}\label{5.2}
{\rm In Proposition \ref{5.1}, $Q$ is assumed to be non-trivial. If we assume further that $M$ induces a Morita equivalence, then the functor $M\otimes_{kHc}-$ induces a bijection
$$\Psi_1:\Proj(kHc)\to \Proj(kGb).$$
For any $(g,h)\in G\times H$, the maps $i_g\Psi_1i_h$ and $\Psi_1$ are the same because both $i_g:\Proj(kGb)\to \Proj(kGb)$ and $i_h:\Proj(kHc)\to \Proj(kHc)$ are identity maps.}
\end{remark}

\medskip\noindent{\it Proof of Corollary \ref{functorialequ}.} This follows from Lemma \ref{lem:criterion}, Proposition \ref{5.1} and Remark \ref{5.2}. $\hfill\square$

\bigskip\noindent\textbf{Acknowledgements.} I wish to thank Professor Yuanyang Zhou for numerous insightful suggestions and helpful discussions which inspired me to consider the topic of Theorem \ref{main}. I am grateful to Doctor Deniz Y{\i}lmaz for kindly explaining my questions on the proof of \cite[Lemma 3.2]{Yilmaz}. I also thank the referee for careful reading and helpful suggestions.

\end{document}